\input amssym.tex

\mag=1200 
\hsize=130mm  \vsize=180mm  \voffset 5mm
\lineskiplimit=8pt \lineskip=8pt
\tolerance=10000 \pretolerance=1000 \parindent=0mm \raggedright

 \font\uh =cmr10 at  18pt

\font\tengo=eufm10

\font\sevengo=eufm7

\font\fivego=eufm5

\font\tenbb=msbm7 at 10pt

\font\sevenbb=msbm7   
 
\font\fivebb=msbm5

%%%%%% Fontes gothique et "blackboard" %%%%%%

\newfam\gofam  \textfont\gofam=\tengo
\scriptfont\gofam=\sevengo   \scriptscriptfont\gofam=\fivego
\def\go{\fam\gofam\tengo}

\newfam\bbfam  \textfont\bbfam=\tenbb
\scriptfont\bbfam=\sevenbb   \scriptscriptfont\bbfam=\fivebb
\def\bb{\fam\bbfam\tenbb}

%%%%%%%%%%%%%%%%%%%%%%%%%%%%%%%%%%%%%%%%%%%%%%%%%%%%%

\def\ind{\hskip 1em\relax}

\def\dim{\mathop{\rm dim}\nolimits}

\def\Ext{\mathop{\rm Ext}\nolimits}
 \def\C{{\bb C}}  \def\Z{{\bb Z}} 
\def\cl{\centerline}
\def\{{\lbrace}  
\def\}{\rbrace}   
\def\({\langle}  
\def\){\rangle}
\def\[{\lbrack} 
\def\]{\rbrack}

\def\arrow{\rightarrow}

\def\.{\bullet}
\def\bs{\bigskip}

\def\ds{\mathop{\oplus}}

%
%%%%%%%%%%%%%%%%%%%%%%%%%%%%%%%%%%%%%%%%%%%%%%%%%%%%%%%%%%%%%%%%%%%%%
%
\cl{\uh A \ naive \ question \ about \ quantum \ groups} \bs \bs 
Let $G$ be a connected semisimple Lie group with finite center~; consider, 
using standard notation, its category $\cal O \subset \go g$-mod of 
BGG-modules, its category $\cal H$ of Harish-Chandra modules, its (complex) 
flag variety $G_{\C}/B$, its compact symmetric space 
$G_{\hbox{\sevenrm c}}/K$ --- and recall the following theorems. \bs
\vbox{{\parindent=7mm
\item{(1)} {\bf Theorem} (BGG). For any simple finite dimensional object $V$ of 
$\cal O$ there is a graded algebra isomorphism
$$\Ext^{\.}_{\cal O}(V,V) \simeq H^{\.}(G_{\C}/B,\C).$$ \par}}  
{\parindent=7mm
\item{(2)} {\bf Theorem} (\'E. Cartan, Casselman). For any simple finite 
dimensional object $V$ of $\cal H$ there is a graded algebra isomorphism
$$\Ext^{\.}_{\cal H}(V,V) \simeq H^{\.}(G_{\hbox{\sevenrm c}}/K,\C).$$ \par}
I think of these statements as being some kind of cohomological Schur Lemmas, 
whence the following definition. \bs
\vbox{{\parindent=7mm
\item{(3)} {\bf Definition}. Let $X$ be a topological space and $\cal A$ be a 
$\C$-category [see Bass [1] p. 57] equipped with a functor 
$F: \cal A \arrow \C$-mod. Then $\cal A$ is a {\bf Schur category} over $X$ if 
$$\left.\matrix{V \in {\cal A} \cr \cr
V \hbox{ simple } \cr \cr
\dim FV < \infty \cr}
\right\} \ \Longrightarrow \ {\Ext}^{\.}(V,V) \simeq H^{\.}(X,\C)$$ 
[isomorphism of graded algebras]. \par}} \bs
In this terminology Theorems (1) and (2) take the respective forms 
``$\cal O$ is a Schur category over $G_{\C}/B$" and ``$\cal H$ is a Schur 
category over $G_{\hbox{\sevenrm c}}/K$". \bs
\ind The purpose of these few lines is to present a conjectural quantum 
analog of Theorem (1). To this end I proceed in two steps. First I define a 
category, denoted ${\cal O}({\go g},h,f)$, which is supposed to be a quantum 
analog of the category $\cal O$ [or more precisely of the category $\cal O$ 
``with weights in the root lattice"]~; then I conjecture that 
${\cal O}({\go g},h,f)$ is a Schur category over the flag variety of $\go g$. 
The category ${\cal O}({\go g},h,f)$ will appear as a subcategory of a 
certain category ${\cal C}({\go g},h,f)$, which is itself a quantum analog 
of $(\go g,h)$-mod [or more precisely of the category of $(\go g,h)$-modules 
with weights in the root lattice]. Here are the details. \bs
\vbox{\ind Let

\cl{$\go g$ be a semismple Lie algebra,}  \bs
\cl{$\alpha_{1},..., \alpha_{r}$ a basis of simple roots,}  \bs
\cl{$(a_{ij})$ the Cartan matrix ({\it i.e.}  
$a_{ij}= 2 (\alpha_i | \alpha_j)/(\alpha_i | \alpha_i)$),}  \bs
\cl{$h$ a complex number,} \bs
\cl{$f = (f_1,...,f_r)$ a list of functions $f_i: \Z^r \rightarrow \C$.}}  \bs 
[It might help the reader to know before hand that the classical case will be 
obtained by putting $f_j(n) = \sum_i \ a_{ij} \ n_i \,$.] \bs
\ind Here starts the {\bf definition of the category} ${\cal C}({\go g},h,f)$. \bs
\ind An object $V$ of ${\cal C}({\go g},h,f)$ is a direct sum 
$$V = \ds_{n \in {\Z}^r} V(n)$$
of vector spaces equipped with endomorphisms $x_{i}$, $y_{i}$ 
$(1 \leq i \leq r)$ satisfying 
$$x_{i}V(n) \subset V(n+e_i),$$
$$y_{i}V(n) \subset V(n-e_i),$$
$$[x_{i},y_{j}] v = \delta_{ij} \ f_{j}(n) v \quad \hbox{for} \quad v \in 
V(n),$$ 
where $(e_i)$ is the canonical basis of ${\Z}^r$, and the {\bf quantum Serre 
relations}, which putting 
$$b(i,j) = 1 - a_{ij},$$ 
$$q(i) = \exp \left( (\alpha_{i} | \alpha_{i}) \ {h \over 2} \right),$$
$$z_{i} = x_{i} \ \forall \ i \ \ \hbox{or} \ \ z_{i} = y_{i}  \ \forall \ i,
$$ 
take the form
$$\sum_{k=0}^{b(i,j)} (-1)^{k} \pmatrix{b(i,j) \cr k}_{q(i)} 
z_{i}^{k} \ z_{j} \ z_{i}^{b(i,j)-k} = 0 \quad \forall \ i \not = j.$$
[The classical case is of course given by $h = 0$.] \bs
\ind The morphisms are the obvious ones. [Here ends the definition of the 
category ${\cal C}({\go g},h,f)$.] \bs
{\parindent=7mm
\item{(4)} {\bf Definition of the category} ${\cal O}({\go g},h,f)$. Let 
$U_{h}({\go n})$ be the algebra generated by the $x_i$ subject to the quantum 
Serre relations. Then ${\cal O}({\go g},h,f)$ is the full subcategory of 
${\cal C}({\go g},h,f)$ whose objects are \hbox{$U_h(\go n)$-finite} and of 
finite length. \par} \bs
\ind If $\cal C$ is $\C$-category and $\cal B$ a full sub-$\C$-category, say 
that $\cal B$  is {\bf Ext-full} in $\cal C$ if for all $V, W \in {\cal B}$ 
the natural morphism 
$$\Ext_{\cal B}^\.(V,W) \rightarrow \Ext_{\cal C}^\.(V,W)$$
is an isomorphism. \bs
\vbox{{\parindent=7mm
\item{(5)} {\bf Conjectures}. \medskip
\itemitem{(a)} The categories ${\cal O}({\go g},h,f)$ and 
${\cal C}({\go g},h,f)$ are Schur categories \hbox{[see (3)]} over the flag 
variety of $\go g$, \medskip
\itemitem{(b)} the inclusion 
${\cal O}({\go g},h,f) \subset {\cal C}({\go g},h,f)$ is Ext-full. \par}} \bs
In the classical case [{\it i.e.} \hbox{$h=0$,} 
\hbox{$f_j(n) = \sum_i a_{ij} \, n_i$]} (a) is due to BGG [see Theorem~(1)]. 
Fuser checked the conjecture for 
\hbox{${\go g = sl}(2,\C)$.} --- Let $\cal C$ be either 
${\cal C}({\go g},h,f)$ or ${\cal O}({\go g},h,f)$ and
\hbox{$\{V_i \ | \ i \in I\}$} a system of representatives of the simple 
objects in $\cal C$.   \bs
{\parindent=7mm
\item{(6)} {\bf Conjecture}. The vector space 
$\ds_{p,i,j} \Ext_{\cal C}^{p}(V_i,V_j)$ is a [nonunital] Koszul algebra.
\par} \bs
This conjecture has been proved for ${\go sl}(2,\C)$ by Fuser and for the 
classical category $\cal O$ by Beilinson, Ginzburg and Soergel (see [2]). 
\bs \cl{* \ * \ * } \bs 
{\par \parindent=6mm 
\item{[1]} Bass H., {\bf Algebraic K-theory}, Benjamin, New York 1968. \bs
\item{[2]} Beilinson A., Ginzburg V., Soergel W., Koszul duality patterns in 
representation theory, {\it J. Am. Math. Soc.} {\bf 9} No.2  (1996) 473-527. \par}

\bigskip

This text and others are available at http://www.iecn.u-nancy.fr/$\sim$gaillard
\bye